\documentclass[11pt]{amsart}
\usepackage{graphicx}
\newtheorem{thm}{Theorem}

\newcommand{\R}{\mathbb{R}}

\newcommand{\E}{\mathbb{E}}
\newcommand{\F}{\mathcal{F}}
\newcommand{\G}{\mathcal{G}}
\newcommand{\V}{\mathcal{V}}
\newcommand{\M}{\mathcal{M}}
\newcommand{\N}{\mathbb{N}}
\newcommand{\B}{\mathcal{B}}
\renewcommand{\L}{\mathcal{L}}
\newcommand{\K}{\mathcal{K}}
\author[1]{Andrew Frohmader}
\author[2]{Hans Volkmer}

\begin{document}
\title[]{1-Wasserstein Distance on the Standard Simplex}

\maketitle

\begin{abstract}
    Wasserstein distances provide a metric on a space of probability measures. We consider the space $\Omega$ of all probability measures on the finite set $\chi = \{1, \dots ,n\}$ where $n$ is a positive integer. 1-Wasserstein distance, $W_1(\mu,\nu)$ is a function from $\Omega \times \Omega$ to $[0,\infty)$. This paper derives closed form expressions for the First and Second moment of $W_1$ on $\Omega \times \Omega$ assuming a uniform distribution on $\Omega \times \Omega$.
\end{abstract}

\section{Introduction}

Wasserstein distances provide a natural metric on a space of probability measures. Intuitively, they measure the minimum amount of work required to transform one distribution into another. They have found application in numerous fields including computer vision \cite{Rubner2000}, statistics \cite{Panaretos19}, dynamical systems \cite{villani2008optimal}, and many others.  

Let $(\chi, d)$  be a Polish space and $\mu$ and $\nu$ be two probability measures on $\chi$. The \textit{p}-Wasserstein distance between $\mu$ and $\nu$ is defined by:

\begin{equation*}
    W_p(\mu,\nu) = \bigg ( \inf_{\pi \in \Pi(\mu, \nu)} \int_{\chi \times \chi} d(x,y)^p d\pi(x,y) \bigg )^{1/p},
\end{equation*}

\noindent where $\Pi(\mu, \nu)$ is the set of probability measures $\pi$ on $\chi \times \chi$ such that $\pi(A \times \chi) = \mu(A)$ and $\pi(\chi \times B) = \nu(B)$ for all measurable $A,B \subseteq \chi$. The measure $\pi$ with marginals $\mu$ and $\nu$ on the two components of $\chi \times \chi$ is called a \textit{coupling} of $\mu$ and $\nu$.

In general,the Wasserstein distances do not admit closed form expressions. One important exception is the case where $\chi = \R$, $d(x,y) = |x-y|$, and $p = 1$. In this case, we have the explicit formula \cite{Panaretos19}:

\begin{equation*}
    W_1(\mu, \nu) = \int_\R |F_\mu(t) - F_\nu(t)|dt,
\end{equation*}

\noindent where $F_\mu(t)$ and $F_\nu(t)$ are the cumulative distribution functions of $\mu$ and $\nu$ respectively. 

Here we specialize further, considering only $\chi = \{1, \dots ,n\}$ where $n$ is a positive integer. Denote this finite set $[n]$. Let $\mu$ and $\nu$ be two probability measures on $[n]$. Then $W_1(\mu, \nu)$ reduces to $l_1$ distance between finite vectors:

\begin{equation} \label{eq_W1_closed}
    W_1(\mu, \nu) = \sum_{i=1}^n |F_\mu(i) - F_\nu(i)| = ||F_\mu - F_\nu||_1.
\end{equation}

There are numerous applications of Wasserstein distances on finite sets. In \cite{Bourn2019}  Bourn and Willenbring, use $W_1$, which they refer to as \textit{Earth Movers Distance} or EMD following computer science convention, to compare grade distributions - distributions on the finite set $\{A, A-, B+, \dots , F\}$. In that context, Bourn and Willenbring considered the space of all probability distributions on $[n]$, denoted $\mathcal{P}_n$ - the \textit{standard simplex} or \textit{probability simplex}. This space embeds naturally in $\R^n$, inherits Lebesgue measure, and has finite volume. A uniform probability measure is obtained by normalizing Lebesgue measure on $\mathcal{P}_n$ so that total mass of $\mathcal{P}_n$ is one. Similarly, $\mathcal{P}_n \times \mathcal{P}_n$ can be given a uniform product probability measure. Let $\mathcal{P}_n \times \mathcal{P}_n = \Omega_n$ with normalized Lebesgue measure be our probability space. Define the random variable $X_n(\mu, \nu) = W_1(\mu, \nu)$ on $\Omega_n$. What can we say about $X_n$?

The contribution of \cite{Bourn2019} to this question is the derivation of a recurrence $\M_{p,q}$ such that $\E(X_n) = \M_{n,n}$. The recurrence is:

\begin{equation} \label{eq_bourn}
    \M_{p,q} = \frac{(p-1)\M_{p-1,q}+(q-1)\M_{p,q-1}+|p-q|}{p+q-1},
\end{equation}

\noindent with $\M_{p,q} = 0$ if either $p$ or $q$ is not positive. In this paper, we find a closed form for the first and second moment of $X_n$. Additionally, we ``solve'' the recurrence provided by \cite{Bourn2019} to verify consistency of results. This provides an interesting link between the discrete combinatorial approach of \cite{Bourn2019} and our calculus based approach.

The remainder of this paper is organized as follows. Section \ref{sec_discussion} presents our results and provides some discussion. Section \ref{sec_approach} gives an overview of our approach to the problem. Section \ref{sec_first} proves the closed form for the first moment. Section \ref{sec_second} proves the closed form for the second moment. Section \ref{sec_recurrence} solves the recurrence \eqref{eq_bourn} of Bourn and Willenbring.

\section{Main Result and Discussion}\label{sec_discussion}

The main results of this paper are the following two theorems. 

\begin{thm} \label{thm_first}
    The first moment of $X_n$ is given by
    \begin{equation*}
        \E(X_n) = \frac{2^{2n-3}(n-1)}{(2n-1)!}(n-1)!^2.
    \end{equation*}
\end{thm}

\begin{thm} \label{thm_second}
    The second moment of $X_n$ is given by
    \begin{equation*}
        \E(X_n^2) = \frac{(n-1)(7n-4)}{30n}.
    \end{equation*}
\end{thm}

\subsection{Unit Normalized $W_1$}

\cite{Bourn2019} defines \textit{unit normalized EMD} on $\mathcal{P}_n$ as EMD scaled by $\frac{1}{n-1}$ - this is $W_1(\mu, \nu)\frac{1}{n-1}$. We will refer to this as \textit{unit normalized} $W_1$ or $\widetilde{W_1}$. This makes $P_n$ into a metric space with diameter one. Intuitively, $W_1$ places the finite set $[n]$ on the first $n$ positive integers of the real number line, while $\widetilde{W_1}$ partitions the unit interval $[0,1]$ with $n$ equally spaced points. $W_1$ has a constant distance of one between adjacent points. $\widetilde{W_1}$ has constant maximum distance (diameter) of one. If we define $\widetilde{X_n} = \frac{X_n}{n-1}$, we can write closed forms for the first and second moment of unit normalized $W_1$ on $\Omega$:

$$\E(\widetilde{X_n}) =  \frac{2^{2n-3}(n-1)!^2}{(2n-1)!}.$$
$$\E(\widetilde{X_n^2}) = \frac{(7n-4)}{30n(n-1)}.$$

\subsection{Asymptotic Behavior}
We briefly consider the asymptotic behavior of the first and second moment. The table below presents approximate numeric values for the first and second moment, and the variance, for $X_n$ and $\widetilde{X_n}$. 

\begin{figure}
    \centering
    \includegraphics[width=10cm]{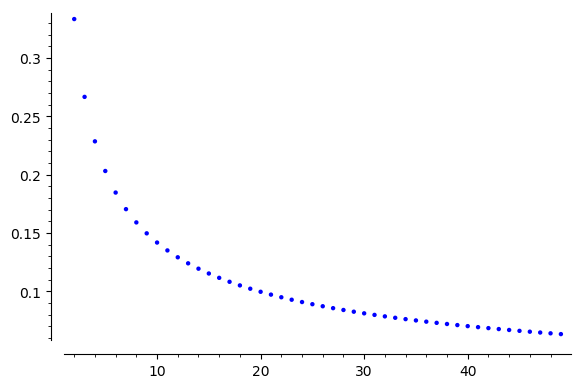}
    \caption{$\E(\widetilde{X_n})$ for n=2 to 50}
    \label{fig:E(X_n)_norm}
\end{figure}

\begin{figure}
    \centering
    \includegraphics[width=10cm]{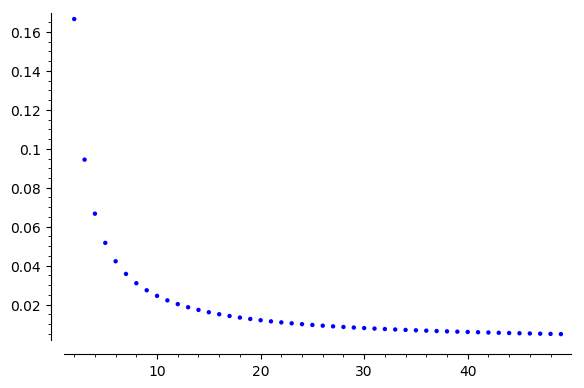}
    \caption{$\E(\widetilde{X_n^2})$ for n=2 to 50}
    \label{fig:E(X_n^2)_norm}
\end{figure}

\renewcommand{\arraystretch}{1.4}
\begin{center}
\begin{tabular}{|l||l|l|l|||l|l|l|} \hline
$n$ & $\E(X_n)$ & $\E(X_n^2)$ & $\text{Var}(X_n)$ & $\E(\widetilde{X_n})$ & $\E(\widetilde{X_n^2})$ & $\text{Var}(\widetilde{X_n})$ \\ \hline
$2$ & $0.3333$ & $0.1667$ & $0.05556$ & $0.3333$ & $0.1667$ & $0.05556$ \\ \hline
$3$ & $0.5333$ & $0.3778$ & $0.09333$ & $0.2667$ & $0.09444$ & $0.02333$ \\ \hline
$4$ & $0.6857$ & $0.6000$ & $0.1298$ & $0.2286$ & $0.06667$ & $0.01442$ \\ \hline
$5$ & $0.8127$ & $0.8267$ & $0.1662$ & $0.2032$ & $0.05167$ & $0.01039$ \\ \hline
$6$ & $0.9235$ & $1.056$ & $0.2027$ & $0.1847$ & $0.04222$ & $0.008107$ \\ \hline
$7$ & $1.023$ & $1.286$ & $0.2392$ & $0.1705$ & $0.03571$ & $0.006645$ \\ \hline
$8$ & $1.114$ & $1.517$ & $0.2759$ & $0.1591$ & $0.03095$ & $0.005630$ \\ \hline
$9$ & $1.198$ & $1.748$ & $0.3126$ & $0.1498$ & $0.02731$ & $0.004884$ \\ \hline
$10$ & $1.277$ & $1.980$ & $0.3493$ & $0.1419$ & $0.02444$ & $0.004313$ \\ \hline
\end{tabular}
\end{center}

Now we consider the behavior of $\E(X_n)$ and $E(\widetilde{X_n})$ as $n \to \infty$. As discussed above:
\begin{equation*}
    \E(X_n) = \frac{2^{2n-3}(n-1)}{(2n-1)!}(n-1)!^2.
\end{equation*}
We can rewrite this in the form:

\begin{equation*}
    \E(X_n) = 4^{n-1}\frac{(n-1)}{n}\frac{1}{\binom{2n}{n}}.
\end{equation*}

It is known that \cite{cameron_1994}

\begin{equation*}
    \binom{2n}{n} \sim \frac{4^n}{\sqrt{\pi n}}.
\end{equation*} 

The symbol $\sim$ means the quotient of the left and right side converges to 1 as $n \to \infty$. Therefore, 
\begin{equation*}
    \E(X_n) \sim \frac{\sqrt{\pi n}}{4}.
\end{equation*} 

Dividing through by $n-1$,
$$ \E(\widetilde{X_n}) \sim \frac{1}{4}\sqrt{\frac{\pi}{n}}.$$

We also have:
$$\E(X_n^2) = \frac{(n-1)(7n-4)}{30n}$$
$$\E(\widetilde{X_n^2}) = \frac{(7n-4)}{30n(n-1)}.$$

Therefore, for the variance we obtain:
$$ \text{Var}(X_n) = \E(X_n^2) - \E(X_n)^2 \sim \bigg (\frac{7}{30} - \frac{\pi}{16} \bigg )n$$

$$ \text{Var}(\widetilde{X_n}) = \E(\widetilde{X_n^2}) - \E(\widetilde{X_n})^2 \sim \bigg (\frac{7}{30} - \frac{\pi}{16} \bigg )\frac{1}{n}.$$

\subsection{Probability Measures for $X_n$}
It would be nice to have explicit probability measures for $X_n$. Toward this end, we computed the probability density functions for $X_2$, $X_3$ as:
$$ p_{X_2}(t) = 2- 2t$$
\[
  p_{X_3}(t) =
  \begin{cases}
       \frac{19}{2}t^4 - \frac{14}{3}t^3+4t^2 & \text{if $0 \leq t \leq 1$} \\
        -\frac{1}{12}t^4+\frac{2}{3}t^3-2t^2+\frac{8}{3}t-\frac{1}{3} & \text{if $1 < t \leq 2$}
  \end{cases}
\]

Further work could be done to determine probability measures for $n>3$.

\begin{figure}
    \centering
    \includegraphics[width=10cm]{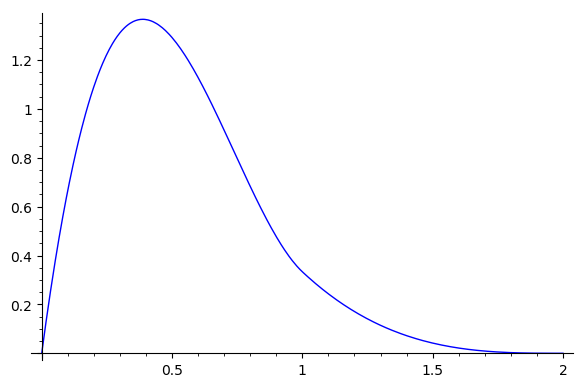}
    \caption{Probability Density Function $p_{X_3}(t)$}
    \label{fig:p_X_3}
\end{figure}

\section{Approach}\label{sec_approach}
We have the probability space $(\Omega_n, \B(\Omega_n), \phi)$ where $\B(\Omega_n)$ denotes the Borel $\sigma$-algebra on $\Omega_n$ and $\phi$ is Lebesgue measure normalized such that $\phi(\Omega_n) = 1$. Our random variable is $X_n(\mu, \nu) = W_1(\mu,\nu)$. It is easy to write the first and second moment of $X_n$ as integrals over $\Omega_n$. The first moment is given by:

\begin{equation*}
    \E(X_n) = \int_{\Omega_n} X_n(\mu, \nu) d\phi.
\end{equation*}

The second moment is given by:

\begin{equation*}
    \E(X_n^2) = \int_{\Omega_n} X_n^2(\mu, \nu) d\phi.
\end{equation*}

For computation, we work with the closed form (\ref{eq_W1_closed}). Let $\mathcal{C}_n$ be the space of cumulative distributions on $[n]$. Consider the probability space $(U_n,\B(U_n), \phi)$, where $U_n = \mathcal{C}_n \times \mathcal{C}_n$,  $\B(U_n)$ denotes the Borel $\sigma$-algebra on $U_n$, and $\phi$ is Lebesgue measure normalized such that $\phi(U_n) = 1$. Define the random variable $Y_n(F_\mu, F_\nu) = ||F_\mu - F_\nu||_1$ on $U_n$. The transformation $T: \mathcal{P}_n \to \mathcal{C}_n$ is volume preserving. By this fact and \eqref{eq_W1_closed}, we conclude $\E(X_n) = \E(Y_n)$ and $\E(X_n^2) = \E(Y_n^2)$. Section \ref{sec_first} proves Theorem \ref{thm_first} from the integral:

\begin{equation*}
    \E(Y_n) = \int_{U_n} Y_n(F_\mu, F_\nu) d\phi.
\end{equation*}

Section \ref{sec_second} proves Theorem \ref{thm_second} from the integral:

\begin{equation*}
    \E(Y_n^2) = \int_{U_n} Y_n^2(F_\mu, F_\nu) d\phi.
\end{equation*}

To find the closed forms presented in Theorem \ref{thm_first} and Theorem \ref{thm_second}, we are led to the evaluation of certain sums for which \cite{A=B} is a general reference. 

\section{Proof of Theorem \ref{thm_first}}\label{sec_first}

The space of all cumulative distributions on $[n]$ is given by 
$$\mathcal{C}_n = \{(x_1,...,x_n) \, | \, 0 \leq x_1 \leq ... \leq x_n = 1\}.$$ 

$\mathcal{C}_n$ is a simplex contained in the hyperplane $x_n = 1$ and thus embeds naturally in $\mathbb{R}^{n-1}$. Henceforth, we will think of $\mathcal{C}_n$ as embedded in $\mathbb{R}^{n-1}$. For $t\ge 0$, $n\in\N,$ and $n\ge 2$, let $\mathcal C_n(t)$ be the set:
\[ \mathcal C_n(t)=\{(x_1,x_2,\dots,x_{n-1})\in\R^{n-1}: 0\le x_1\le x_2\le \dots \le x_{n-1}\le t\} .\]

Let $s,t\ge 0$, $\textbf{x} \in \mathcal{C}_n(t)$, $\textbf{y} \in \mathcal{C}_n(s)$ and set

\begin{eqnarray*}
 \V_n(s,t)&=&\int_{\mathcal C_n(s)}\int_{\mathcal C_n(t)} 1\, d \textbf{x} \, d \textbf{y}=\frac{s^{n-1}t^{n-1}}{(n-1)!^2},\\
 \F_n(s,t)&=&\int_{\mathcal C_n(s)}\int_{\mathcal C_n(t)} \|\textbf{x} - \textbf{y}\|_1 d \textbf{x} \, d \textbf{y}.\\
\end{eqnarray*}
Then,

$$ \E(Y_n) = \frac{\F_n(1,1)}{\V_n(1,1)}.$$

If we set $\F_1(s,t)=0$, then the functions $\F_n$ are determined recursively by
\[
\F_n(s,t)= \int_0^t\int_0^s (\F_{n-1}(x,y)+|x-y|\V_{n-1}(x,y))\,dx\,dy.
\]
Here $x = x_{n-1}$ and $y=y_{n-1}$ are the last components of $\mathcal{C}_n(t)$ and $\mathcal{C}_n(s)$ respectively. Fubini's theorem is relied on to justify evaluating $\F_n(s,t)$ as an iterated integral. 

This recursion formula is awkward to use because of the appearance of the absolute value $|x-y|$.
We remove this absolute value as follows.
The function $\F_n(s,t)$ is symmetric in $(s,t)$.
Therefore, it is sufficient to consider $s\ge t$ which we will always assume.
Introduce the integral operator
\[ (I f)(s,t)=2\int_0^t\int_y^t f(x,y)\,dx\,dy+\int_0^t\int_t^s f(x,y)\,dx\,dy .\]
Then, for $0\le t\le s$,
\begin{equation}\label{recF}
\F_n(s,t)= I\left(\F_{n-1}(x,y)+(x-y)\V_{n-1}(x,y)\right).
\end{equation}
This is a recursion formula for $\F_n(s,t)$ in the region $0\le t\le s$.
We note that, for $i,j=0,1,2,\dots$,
\begin{equation}\label{I}
 I(x^iy^j)=\frac{s^{i+1}t^{j+1}}{(i+1)(j+1)}+\frac{(i-j)t^{i+j+2}}{(i+1)(j+1)(i+j+2)} .
\end{equation}
It follows from \eqref{recF} and \eqref{I} that $\F_n(s,t)$ is a polynomial in $s,t$ of the form
\begin{equation}\label{F}
\F_n(s,t)=\sum_{k=0}^{n} f_{n,k} t^{n+k-1}s^{n-k} .
\end{equation}
For example, we get
\begin{eqnarray*}
\F_2(s,t)&=&\frac{ts^2}{2}-\frac{t^2s}{2}+\frac{t^3}{3},\\
\F_3(s,t)&=&\frac{t^2s^3}{4}-\frac{t^3s^2}{4}+\frac{t^4s}{12}+\frac{t^5}{20},\\
\F_4(s,t)&=&\frac{t^3s^4}{24}-\frac{t^4s^3}{24}+\frac {t^5s^2}{120}+\frac {t^{6}s}{120}+\frac {t^7}{420},\\
\F_5(s,t)&=&\frac {t^4s^5}{288}-\frac{t^5s^4}{288}+\frac {t^6s^3}{2160}+\frac {{t}^{7}{s}^{
2}}{1680}+\frac {t^8s}{3360}+\frac {t^9}{18144}.
\end{eqnarray*}

Using \eqref{recF} and \eqref{I}, we find
\[ f_{n,0}=\frac{f_{n-1,0}}{n(n-1)}+\frac{1}{n(n-1)(n-2)!^2} \quad \text{for $n\ge 2$}.\]
We set
\[ \tilde f_{n,k}=f_{n,k}(n+k-1)!(n-k)! \]
Then
\[ \tilde f_{n,0}-\tilde f_{n-1,0} =n-1.\]
Since $f_{1,0}=0$ we find
$\tilde f_{n,0}=\frac12 n(n-1)$ so
\begin{equation}\label{f1}
f_{n,0}=\frac{n(n-1)}{2(n-1)!n!} =\frac{1}{2(n-1)!(n-2)!} .
\end{equation}
Similarly, we obtain
\begin{equation}\label{f2}
f_{n,1}=-f_{n,0}=-\frac{1}{2(n-1)!(n-2)!} .
\end{equation}
Now \eqref{recF} and \eqref{I} give
\[ f_{n,k}=\frac{f_{n-1,k}}{(n+k-1)(n-k)} \quad\text{for $2\le k<n$.}\]
Therefore,
\[ \tilde f_{n,k}=\tilde f_{n-1,k} \quad\text{for $2\le k<n$}\]
and
\begin{equation}\label{f3}
f_{n,k}=\frac{f_{k,k} (2k-1)!}{(n+k-1)!(n-k)!} \quad\text{for $2\le k\le n$.}
\end{equation}
Using again \eqref{recF} and \eqref{I} we obtain
\begin{equation}\label{f4}
f_{n,n}=\sum_{k=0}^{n-1} \frac{f_{n-1,k}(1-2k)}{(n+k-1)(n-k)(2n-1)} +\frac{2}{n!(n-2)!(2n-1)}.
\end{equation}
In \eqref{f4} we substitute \eqref{f1}, \eqref{f2}, \eqref{f3} and find
\begin{equation}\label{f5}
f_{n,n}=\sum_{k=2}^{n-1} \frac{f_{k,k} (2k-1)!(1-2k)}{(n+k-1)!(n-k)!(2n-1)} +\frac{1}{(n-1)!(n-2)!(2n-1)} .
\end{equation}
This is a recursion formula for the ``diagonal sequence'' $f_{k,k}$. With the help of a computer we guess the solution
\begin{equation}\label{f6}
 f_{k,k}=\frac{(k-1)k}{(2k-1)!} .
 \end{equation}
We prove \eqref{f6} by induction on $k$. The formula is true for $k=1$. Suppose \eqref{f6} holds for
all $k$ less than a given $n\ge 2$. Then we substitute $f_{k,k}$ given by \eqref{f6} on the right-hand side of \eqref{f5}.
Since (see Appendix \ref{app_binomial_sums})
\[ \sum_{k=2}^n \frac{k(k-1)(2k-1)}{(n+k-1)!(n-k)!}=2A_{n,3}-3A_{n,2}+A_{n,1}=\frac{1}{(n-1)!(n-2)!} \]
we obtain \eqref{f6} for $k=n$. This completes the inductive proof of \eqref{f6}.

Now \eqref{f3}, \eqref{f6} give
\begin{equation}\label{f7}
f_{n,k}=\frac{(k-1)k}{(n+k-1)!(n-k)!} \quad \text{for $2\le k\le n$} .
\end{equation}

In particular, using again Appendix \ref{app_binomial_sums}, we have
\begin{eqnarray*}
\F_n(1,1)&=&\sum_{k=0}^n f_{n,k}=\sum_{k=2}^n f_{n,k}\\
&=&\sum_{k=2}^n \frac{(k-1)k}{(n+k-1)!(n-k)!}\\
&=& A_{n,2}-A_{n,1}\\
&=&\frac{2^{2n-3} (n-1)}{(2n-1)!} .
\end{eqnarray*}

Thus,

$$ \E(X_n) = \E(Y_n) = \frac{\F_n(1,1)}{\V_n(1,1)} = \frac{2^{2n-3} (n-1)}{(2n-1)!}(n-1)!^2.$$

\section{Proof of Theorem \ref{thm_second}}\label{sec_second}

Let $s,t\ge 0$, $\textbf{x} \in \mathcal{C}_n(t)$, $\textbf{y} \in \mathcal{C}_n(s)$ and define
\[\G_n(s,t)=\int_{\mathcal C_n(s)}\int_{\mathcal C_n(t)} \|\textbf{x}-\textbf{y}\|_1^2 d\textbf{x} \, d\textbf{y}.
\]
If we set  $\G_1(s,t)=0$
then the function $\G_n$ is determined recursively by
\[\G_n(s,t)= \int_0^t\int_0^s (\G_{n-1}(x,y)+2|x-y|\F_{n-1}(x,y)+(x-y)^2\V_{n-1}(x,y))\,dx\,dy.
\]
As in the proof of Theorem \ref{thm_first}, $x = x_{n-1}$, $y=y_{n-1}$ and Fubini's theorem justifies evaluating $\G_n(s,t)$ as an iterated integral.

The function $\G_n(s,t)$ is symmetric in $(s,t)$.
Therefore, it is sufficient to consider $s\ge t$.
Then, for $0\le t\le s$,
\begin{equation}\label{recG}
\G_n(s,t)= I\left(\G_{n-1}(x,y)+2(x-y)\F_{n-1}(x,y)+(x-y)^2\V_{n-1}(x,y)\right).
\end{equation}
This is a recursion formula for $\G_n(s,t)$ in the region $0\le t\le s$.
It follows from \eqref{recG} and \eqref{I} that $\G_n(s,t)$ is a polynomial in $s,t$ of the form
\begin{equation}\label{G}
\G_n(s,t)=\sum_{k=-1}^n g_{n,k} t^{n+k}s^{n-k} .
\end{equation}
For example,
\begin{eqnarray*}
\G_2(s,t)&=&\frac{ts^3}{3}-\frac{t^2s^2}{2}+\frac{t^3s}3,\\
\G_3(s,t)&=& \frac{7t^2s^4}{24}-\frac{t^3s^3}{2}+\frac{3t^4s^2}{8}-\frac{2t^5s}{15}+\frac{11t^6}{180},\\
\G_4(s,t)&=&\frac{5t^3s^5}{72}-\frac{t^4s^4}8+\frac{31t^5s^3}{360}-\frac{t^6s^2}{60}-\frac{t^7s}{180}+\frac{t^8}{120},\\
\G_5(s,t)&=&\frac{13 t^4s^6}{1728}-\frac{t^5s^5}{72}+\frac{77 t^6s^4}{8640}-\frac{t^7s^3}{1260}-\frac{11 t^8s^2}{10080}+\frac{t^9s}{2520}+\frac{19t^{10}}{50400}.
\end{eqnarray*}

First we determine $g_{n,k}$ for $k=-1$.
From \eqref{recG} and \eqref{I} we obtain
\[ g_{n,-1}=\frac{g_{n-1,-1}}{(n-1)(n+1)}+\frac{2f_{n-1,0}}{(n-1)(n+1)}+\frac{1}{(n-2)!^2(n-1)(n+1)} .\]
Using \eqref{f1} this yields
\[ g_{n,-1}=\frac{g_{n-1,-1}}{(n-1)(n+1)}+\frac{1}{(n-1)!(n-3)!(n+1)}+\frac{1}{(n-2)!^2(n-1)(n+1)} .\]
Let
\[ g_{n,k}=\frac{\tilde g_{n,k}}{(n-k)!(n+k)!} .\]
Then
\[ \tilde g_{n,-1}-\tilde g_{n-1,-1}=(n-1)^2 n .\]
This gives
\[ \tilde g_{n,-1}=\tilde g_{1,-1}+\sum_{m=2}^n (m-1)^2m =\frac1{12} n(n-1)(3n-2)(n+1).\]
Therefore,
\begin{equation}\label{g1}
g_{n,-1}=\frac{ n(n-1)(3n-2)(n+1)}{12 (n-1)!(n+1)!}=\frac{n(3n-2)}{12(n-2)!n!} .
\end{equation}
From \eqref{recG} and \eqref{I} we obtain the recursion formula
\[ g_{n,0}=\frac{g_{n-1,0}}{n^2} +\frac{2f_{n-1,1}-2f_{n-1,0}}{n^2}-\frac{2}{n^2(n-2)!^2} .\]
This gives
\[ \tilde g_{n,0}-\tilde g_{n-1,0} =-2(n-1)^3\]
and
\begin{equation}\label{g2}
g_{n,0}=-\frac{1}{2(n-2)!^2} .
\end{equation}

From \eqref{recG} and \eqref{I} we obtain the recursion formula
\[ g_{n,1}=\frac{g_{n-1,1}}{(n-1)(n+1)} +\frac{2f_{n-1,2}-2f_{n-1,1}}{(n-1)(n+1)}+\frac{1}{(n-2)!^2(n-1)(n+1)} .\]
This gives
\[ \tilde g_{n,1}-\tilde g_{n-1,1} =4(n-2)+n(n-1)^2\]
and
\begin{equation}\label{g3}
g_{n,1}= \frac{2(n-1)(n-2)}{(n-1)!(n+1)!}+\frac{ n(n-1)(3n-2)(n+1)}{12 (n-1)!(n+1)!} .
\end{equation}
Now let $k\ge 2$. Then \eqref{recG} and \eqref{I} give
\[ g_{n,k}=\frac{g_{n-1,k}}{(n+k)(n-k)}+\frac{2f_{n-1,k+1}-2f_{n-1,k}}{(n+k)(n-k)} .\]
This leads to the recursion formula
\[ \tilde g_{n,k}-\tilde g_{n-1,k}=4k(n-k^2-1)\quad \text{for $2\le k<n$} .\]
This gives
\begin{equation}\label{g4}
\tilde g_{n,k}=\tilde g_{k,k}+2k(k-n)(2k^2-k-n+1).
\end{equation}

From \eqref{recG} and \eqref{I} we find
\begin{equation}\label{g5}
 g_{n,n}=-\sum_{k=-1}^{n-1} \frac{k g_{n-1,k}}{(n-k)(n+k) n}
-\sum_{k=1}^{n-2} \frac{2 kf_{n-1,k+1}}{(n-k)(n+k)n}+\sum_{k=2}^{n-1} \frac{2kf_{n-1,k}}{(n-k)(n+k)n} .
\end{equation}
It should be noted that the term $(x-y)^2 \V_{n-1}(x,y)$ in \eqref{recG} does not contribute
to this formula because it is a symmetric polynomial.
Using \eqref{f7} and Appendix \ref{app_binomial_sums} we evaluate the second and third sum on the right-hand side of \eqref{g5}.
We find
\begin{eqnarray*}
 -\sum_{k=1}^{n-2} \frac{2 kf_{n-1,k+1}}{(n-k)(n+k)n}&=&\frac{4^{n-1}(n+1)}{(2n)!}-\frac{n(n-2)}{n!^2}-
\frac{n+1}{(2n-1)!},\\
\sum_{k=2}^{n-1} \frac{2kf_{n-1,k}}{(n-k)(n+k)n}&=&\frac{4^{n-1}(n+1)}{(2n)!}+\frac{n(n-2)}{n!^2}-\frac1{2(2n-3)!}.
\end{eqnarray*}
Therefore, \eqref{g5} can be written as
\begin{equation}\label{g6}
 g_{n,n}=-\sum_{k=-1}^{n-1} \frac{k g_{n-1,k}}{(n-k)(n+k) n} +\frac{2^{2n-1}(n+1)}{(2n)!}-\frac{2(n^2-n+1)}{(2n-1)!} .
\end{equation}
Using \eqref{g1}, \eqref{g3} we find that
\begin{equation}\label{g7}
-\sum_{k=-1}^1 \frac{k g_{n-1,k}}{(n-k)(n+k) n} =-\frac{2(n-3)(n-2)}{n!^2(n+1)}.
\end{equation}
We also have
\begin{eqnarray}\label{g8}
&&-\sum_{k=2}^{n-1} \frac{2k^2(k-n+1)(2k^2-k-n+2)}{n(n-k)!(n+k)!} \\ &&=\frac{4^{n-1}(n+1)(4n-3)}{(2n)!}+\frac{2(n^2-n+1)}{(2n-1)!}-\frac{4n^3+n^2+9n-12}{(n+1)n!^2}.\nonumber
\end{eqnarray}
We now substitute \eqref{g4}, \eqref{g7}, \eqref{g8} in \eqref{g6} and obtain
\begin{equation}\label{g9}
 g_{n,n}=-\sum_{k=2}^{n-1} \frac{k (2k)! g_{k,k}}{(n-k)!(n+k)! n} +\frac{4^{n-1}(n+1)(4n-1)}{(2n)!} -\frac{4n-1}{n(n-1)!^2}.
\end{equation}
This is a recursion formula for the sequence $g_{n,n}$.
With the help of a computer we guess
\begin{equation}\label{g10}
 g_{k,k}=\frac{(k-2)(4k-1)}{30(2k-3)!}\quad\text{for $k\ge 2$} .
\end{equation}
We prove \eqref{g10} by induction on $k$. This equation is true for $k=2$.
Suppose it is true for all $k$ less than a given $n\ge 3$.
Then we use \eqref{g10} in \eqref{g9}. Using Appendix \ref{app_binomial_sums}, we evaluate the sum and and obtain the correct expression for $g_{n,n}$. The inductive proof is complete.

From \eqref{g4} we get, for $2\le k\le n$,
\begin{equation}\label{g11}
g_{n,k}=\frac{k(2k-1)(2k-2)(k-2)(4k-1)}{15 (n+k)!(n-k)!}+\frac{2k(k-n)(2k^2-k-n+1)}{(n+k)!(n-k)!}.
\end{equation}
In particular, we have
\[ \G_n(1,1)=\sum_{k=-1}^n g_{n,k} =\frac{n(n-1)(7n-4)}{30n!^2} .\]

Thus,

$$ \E(X_n^2) = \E(Y_n^2) = \frac{\G_n(1,1)}{\V_n(1,1)} = \frac{(n-1)(7n-4)}{30n}.$$

\section{Solution to  \cite{Bourn2019} Recurrence} \label{sec_recurrence}

In \cite{Bourn2019} a double sequence $\M_{p,q}$, $p,q\in\N_0$, is defined as the solution of the recursion \eqref{eq_bourn}.
Here we ``solve'' this recursion.

In order to simplify the recursion we set
\[ \L_{p,q}=\frac{(p+q-1)!}{(p-1)!(q-1)!}\M_{p,q},\quad\text{$p,q\in\N$,}\]
and $\L_{p,0}=\L_{0,q}=0$.
Then \eqref{eq_bourn} becomes
\begin{equation}\label{recL}
\L_{p,q}=\L_{p-1,q}+\L_{p,q-1}+|p-q|\binom{p+q-2}{p-1},\quad\text{$p,q\in\N$}.
\end{equation}
Note that $\L_{p,q}$ is a nonnegative integer for all $p,q$.
The matrix $\L_{p,q}$, $p,q=0,\dots,6$, is
\[
\begin{pmatrix}
0 & 0 & 0 & 0 & 0 & 0 & 0\\
0 & 0 & 1 & 3 & 6 & 10 & 15\\
0 & 1 & 2 & 8 & 22 & 47 & 86\\
0 & 3 & 8 & 16 & 48 & 125 & 274\\
0 & 6 & 22 & 48 & 96 & 256 & 642\\
0 & 10 & 47 & 125 & 256 & 512 &1280\\
0 & 15 & 86 & 274 & 642 & 1280 & 2560
\end{pmatrix}
\]

Consider the more general recursion
\begin{equation}\label{recK}
\K_{p,q}=\K_{p-1,q}+\K_{p,q-1}+ a_{p,q} ,\quad p,q\in\N,
\end{equation}
with initial condition $\K_{p,0}=\K_{0,q}=0$, where $\{a_{p,q}\}$ is a given double sequence.
As a special case consider $a_{p,q}=1$ if $p=i$, $q=j$ for some given $i,j\in\N$, and $a_{p,q}=0$ otherwise.
If $p\ge i$, $q\ge j$ the corresponding solution is
\[ \K_{p,q}=\binom{p+q-i-j}{p-i} .\]
Therefore, by superposition, the solution of \eqref{recL} is
\[
\L_{p,q}=\sum_{i=1}^p \sum_{j=1}^q |i-j|\binom{i+j-2}{i-1}\binom{p+q-i-j}{p-i}.
\]
We rewrite this in the slightly nicer form
\begin{equation}\label{L}
\L_{p+1,q+1}= \sum_{i=0}^p \sum_{j=0}^q |i-j|\binom{i+j}{i}\binom{p+q-i-j}{p-i}.
\end{equation}
In a sense we found the solution of \eqref{recL} and so also of \eqref{eq_bourn}. However, the question arises whether
we can express the double sum in \eqref{L} (more) explicitly.
We evaluate the double sum in \eqref{L} when $p=q$. There is some evidence (from the computer)
that this is not possible when $p\ne q$. We could handle some special case like small $p$ or small $|p-q|$.

We now evaluate the double sum in \eqref{L} when $p=q$.
First, we introduce $i+j=m$ as a new summation index. Then
\[ \L_{p+1,p+1}=\sum_{m=0}^{2p} \sum_{i=0}^{p} |m-2i|\binom{m}{i}\binom{2p-m}{p-i} ,\]
where we use the standard convention  that $\binom{n}{k}=0$ unless $0\le k\le n$.
We remove the absolute value by writing
\begin{equation}\label{eq1}
 \L_{p+1,p+1}=2\sum_{m=0}^{2p} \sum_{i=0}^{\lfloor m/2\rfloor} (m-2i)\binom{m}{i}\binom{2p-m}{p-i} .
\end{equation}
The inside sum can be evaluated as a telescoping sum as follows.
Consider first an even $m=2r$. Then
\[ \sum_{i=0}^r 2(r-i)\binom{2r}{i}\binom{2p-2r}{p-i} =\sum_{j=0}^{r} 2j\binom{2r}{r-j}\binom{2s}{s+j},\]
where we set $i=r-j$ and $s=p-r$. Note that $\binom{2r}{r-j}=\binom{2r}{r+j}$.
Now
\[  2j\binom{2r}{r-j}\binom{2s}{s+j}=d_j-d_{j+1},\quad j=0,1,\dots,r,\]
where
\[ d_j=\frac{(r+j)(s+j)}{r+s}\binom{2r}{r-j}\binom{2s}{s+j}.\]
Note that $d_{r+1}=0$. Therefore,
\[ \sum_{i=0}^r 2(r-i)\binom{2r}{i}\binom{2p-2r}{p-i}=d_0-d_{r+1}=d_0=\frac{rs}{r+s}\binom{2r}{r}\binom{2s}{s}.\]
If $m=2r+1$ is odd, we use
\[ (2j+1)\binom{2r+1}{r-j}\binom{2s-1}{s+j}=e_j-e_{j+1},\]
where
\[ e_j=\frac{(r+1+j)(s+j)}{r+s}\binom{2r+1}{r-j}\binom{2s-1}{s+j}.\]
Therefore,
\[ \sum_{i=0}^r (2r+1-2i)\binom{2r+1}{i}\binom{2p-2r-1}{p-i}=e_0=\frac{(r+1)s}{r+s}\binom{2r+1}{r}\binom{2s-1}{s-1}.\]

Substituting these results in \eqref{eq1} we obtain
\begin{equation}\label{eq2}
\L_{p+1,p+1}=\sum_{r=0}^p \frac{2r(p-r)}{p}\binom{2r}{r}\binom{2p-2r}{p-r} +\sum_{r=0}^{p-1} \frac{2(r+1)(p-r)}{p}\binom{2r+1}{r}\binom{2p-2r-1}{p-r-1}.
\end{equation}
The sums appearing in \eqref{eq2} are easy to evaluate. We note that
\[ (1-4x)^{-1/2}=\sum_{n=0}^\infty \binom{2n}{n} x^n\quad\text{for $|x|<\frac14$},\]
so
\[ 2x(1-4x)^{-3/2} =\sum_{n=0}^\infty n\binom{2n}{n} x^n .\]
If we compare coefficients in
\[ \frac{2x}{(1-4x)^{3/2}} \frac{2x}{(1-4x)^{3/2}}=\frac{4x^2}{(1-4x)^3}=\sum_{n=0}^\infty  2^{2n-3} n(n-1) x^n \]
we find the first sum in \eqref{eq2}. Using $\binom{2r+1}{r}=\frac12 \binom{2r+2}{r+1}$ we also find the second sum.
Finally, we obtain
\[ \L_{p+1,p+1}=\frac{2}{p}2^{2p-3}p(p-1)+\frac2p 2^{2p-3}p(p+1). \]
Therefore,
\begin{equation}\label{eq3}
\L_{p+1,p+1}= 2^{2p-1} p,\quad  \M_{p+1,p+1}=2^{2p-1} p \frac{p!^2}{(2p+1)!}.
\end{equation}
Since \cite{Bourn2019} showed that $\E(X_n) = \M(p,p)$, this result agrees with Theorem \ref{thm_first}.

\bigskip

\noindent \textbf{Acknowledgements:} The first author would like to thank Jeb Willenbring and Rebecca Bourn for helpful discussions.

\bibliographystyle{alpha}
\bibliography{references}

\appendix

\section{Binomial Sums} \label{app_binomial_sums}
For $p\in\N_0, n\in\N$ we define
\begin{eqnarray*}
A_{n,p}&=& \sum_{k=1}^n \frac{k^p}{(n+k-1)!(n-k)!},\\
B_{n,p}&=& \sum_{k=1}^n \frac{k^p}{(n+k)!(n-k)!}.
\end{eqnarray*}
These sums can be calculated recursively
from
\begin{eqnarray*}
 A_{n,0}&=&\frac{4^{n-1}}{(2n-1)!},\\
 B_{n,0}&=&\frac{2^{2n-1}}{(2n)!}-\frac{1}{2n!^2} .
\end{eqnarray*}
and
\[ A_{n,p}=nA_{n,p-1}-B_{n-1,p-1}, \quad B_{n,p}=A_{n,p-1}-n B_{n,p-1} .\]
For example, we get
\begin{eqnarray*}
A_{n,1}&=&\frac{2^{2n-3}}{(2n-1)!}+\frac{1}{2(n-1)!^2},\\
A_{n,2}&=&\frac{n2^{2n-3}}{(2n-1)!}+\frac{1}{2(n-1)!^2},\\
A_{n,3}&=&\frac{4^{n-2}(3n-1)}{(2n-1)!}+\frac{1}{2n!(n-1)!},\\
B_{n,1}&=&\frac{1}{2n!(n-1)!},\\
B_{n,2}&=&\frac{2^{2n-3}}{(2n-1)!},\\
B_{n,3}&=&\frac{1}{2(n-1)!^2} .
\end{eqnarray*}.

\end{document}